\documentclass[11pt]{article}
\usepackage[utf8]{inputenc}
\usepackage[
papersize={5.5in, 8.5in},
left=0.35in,
right=0.35in,
top=0.35in,
bottom=0.5in]{geometry}
\usepackage{graphicx}
\usepackage{amsthm}
\usepackage{rotating}
\usepackage{amsfonts,amssymb,amsmath}
\usepackage[usenames]{color}

\usepackage{array}
\usepackage{float}
\usepackage{booktabs}
\usepackage[dvipsnames,usenames,table,xcdraw]{xcolor}

\graphicspath{ {images/} }

\title{\bf
The surface of a sufficiently large sphere has chromatic number at most 7
}
\author{\normalsize {Tomas Sirgedas}} 
\date{\normalsize {Okemos, Michigan, USA, tsirgedas@gmail.com}}

\begin{document}

\maketitle

\begin{abstract}
We present a method to assign, for any radius $r$ greater than about 12.44, one of seven colors to each point in $\mathbb{R}^3$ lying at distance $r$ from the origin, such that no two points at unit distance from each other are assigned the same color. The existence of such a construction contrasts with the recent demonstration that, for any positive value $\varepsilon$, if no two points assigned the same color lie at any distance in $[1,1+\varepsilon]$ (and with certain other restrictions that are also satisfied with our coloring), then eight colors are needed for any finite $r\ge18$, even though seven colors suffice in the plane when $\varepsilon \leq\frac{\sqrt{7}}{2} - 1$.
\end{abstract}

\section{Background}

\subsection{Chromatic number of the plane: upper bounds}

\begin{figure}[!b]
\centering
\includegraphics[scale=0.3]{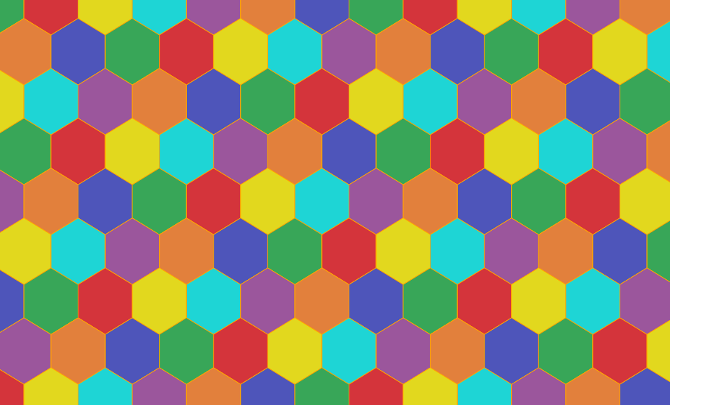}
\caption{The “Isbell tiling” of the plane.}
\label{f1}
\end{figure}

The celebrated “chromatic number of the plane" (CNP) problem, generally termed the Hadwiger-Nelson problem, asks how many colors must be assigned to the points of the real plane so that no pair of points at unit distance is monochromatic. Such an assignment will hereafter be termed a valid coloring. As of this writing, the value of the CNP is known to be $5$, $6$ or $7$. The upper bound of $7$ is easily demonstrated by a "tiling" of the plane (that is, an assignment of colors such that the plane is divided into monochromatic regions, or tiles, each bounded by a continuous closed curve): we can use regular hexagonal tiles of diameter slightly less than $1$, and assign seven colors such that no pair of tiles of the same color is closer than $1$. See Figure \ref{f1}. Hereafter we refer to this as the Isbell tiling, after its discoverer (\cite{soi}, p. 24). More generally, we describe a tiling using $k$ colors as a $k$-tiling.

\subsection{Chromatic number of the plane: conditional lower bounds}

The lower bound on the CNP remained stubbornly stuck at $4$ for decades and has been $5$ since 2018, but higher values have been shown if certain properties of the coloring are stipulated. In particular, in 1999, Thomassen showed \cite{thom} that $7$ is in fact the exact number of colors needed, if we restrict our colorings to those satisfying three conditions, namely:

1) The coloring is a tiling.

2) The tiling is not “Siamese" (our terminology). Let a tile's “annulus of exclusion" (AE) be the set of all points lying at unit distance from any point of the tile. Then a pair of tiles is Siamese if they lie entirely inside the inner boundary of each other's AE. A tiling is Siamese if it includes any Siamese pair of tiles.

\begin{figure}[!b]
\centering
\includegraphics[scale=0.3]{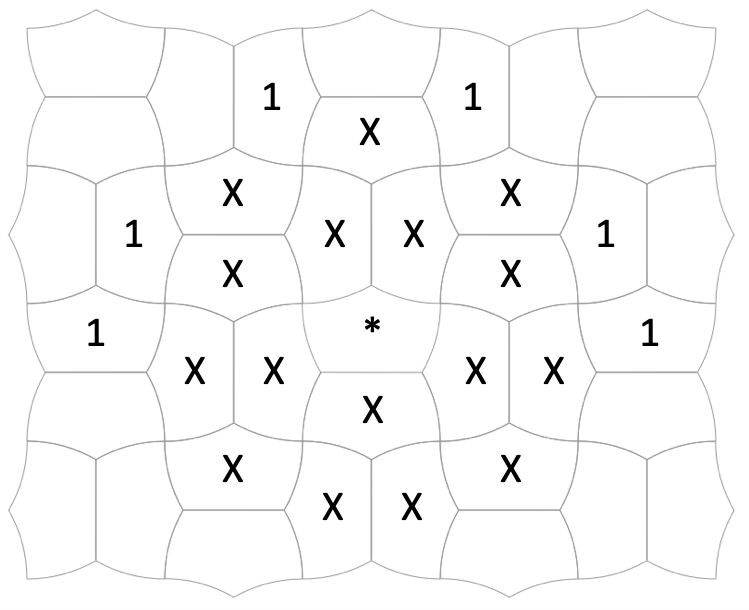}
\caption{An interesting unscalable tiling of the plane. The tiles labeled “X" are too close to the tile labeled “*" to be the same color as it; however, the tiles labeled “1" are exactly 1 from tile “*" and therefore can be the same color as it.}
\label{f2}
\end{figure}

3) The tiling is “scalable". In the Isbell tiling, the ratio of the minimum distance between same-colored tiles to the tile diameter is strictly greater than 1. However, tilings exist in which this is not so, and some of these possess intriguing features. For example, the tiling shown in Figure \ref{f2} employs the tactic of drawing certain inter-tile boundaries as arcs of unit radius, with the result that each tile has diameter exactly 1 and is exactly 1 from several other tiles. Notably, it is less than 1 from only 16 other tiles, whereas in the Isbell tiling the corresponding number is 18. Such observations underpin continuing hope, if not expectation, that a $6$-tiling of the plane exists. Alternatively, it may be possible to enumerate a finite number of cases, along lines resembling the proof of the four-color theorem \cite{aphak}, so as to prove that six colors do not suffice even for unscalable non-Siamese tilings. (It is known \cite{tow} that tilings require at least six colors, whether or not they are Siamese or scalable.) Similarly to “Siamese," we say that a tile is unscalable if it abuts two tiles of the same color (at least one of the two abutments being, necessarily, only a single point) and is thus forced to have a diameter of exactly 1, and that a tiling is unscalable if any component tile is unscalable.

\subsection{Chromatic numbers of other regions}

The CNP question has been adapted to many other spaces, such as $\mathbb{Q}^2$, $\mathbb{R}^3$ and defined subsets of $\mathbb{R}^2$. In this report we are interested in the surface of a sphere in $\mathbb{R}^3$. We define the distance between two points as the Euclidean distance in $\mathbb{R}^3$ rather than the distance on the sphere surface; clearly a sphere colored excluding unit distance with one such convention can be scaled to give one with the other convention, so this choice does not affect what follows.

\begin{figure}[!b]
\centering
\includegraphics[scale=0.35]{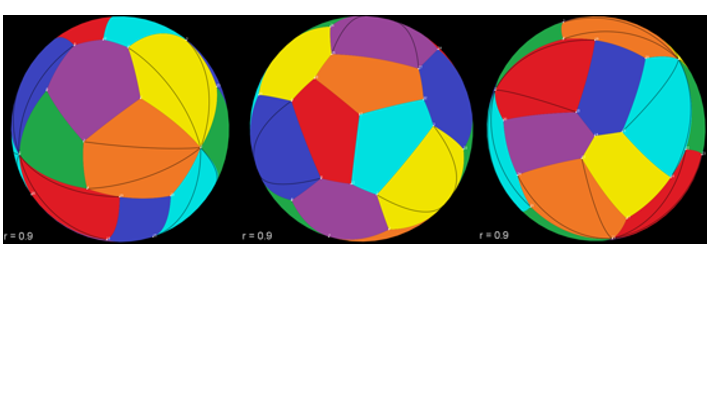}
\caption{A $7$-tiling of a sphere of radius $0.9$. The tiling has $3$-way rotational symmetry about an axis whose poles are at the top and bottom of the projection shown in the left panel. Center: a projection with the North pole in the center. Right: a projection with the South pole in the center. Nine of the tiles share two boundaries, one being only a single point, with tiles of the same color, making the tiling quasi-scalable. Thin black lines indicate tile diameters that must be exactly 1.}
\label{f3}
\end{figure}

For small spheres there are tilings using between $4$ and $6$ colors arising from projection of Platonic solids (or, in the case of $5$ colors, a square pyramid) onto the sphere, and these (see \cite{sim}) remain the best-known upper bounds for the chromatic number of spheres within these ranges of radius. (As an aside, we note that Malen \cite{sim} implies that a sphere of any radius $\sqrt{\frac{1}{3}} < r \leq \frac{\sqrt{3}}{2}$ can be 6-tiled with a projection of a regular dodecahedron, when in fact this only works for $r\ge \sqrt{\frac{3}{8}}$; however, the upper bounds for these ranges remain correct because radii $\sqrt{\frac{1}{3}} < r \leq \sqrt{\frac{3}{8}}$ can be $6$-tiled with a projection of a cube.) We are unaware of any published discussions of tilings of larger spheres, but as an aside from the main theme of this paper we discovered that a sphere of radius between about $0.84$ and $1.005$ can be $7$-tiled with $21$ tiles. See Figure \ref{f3}.

The tiling in Figure \ref{f3} introduces a feature that will crop up extensively below. Like the one in Figure \ref{f2} and the octahedral $4$-tiling of a sphere of radius $\sqrt{\frac{1}{2}}$, it is unscalable, because some tiles have two neighbors of the same color. However, unlike the octahedron case, the radius of the sphere can still vary within a certain range. This is because different radii lead to tiles with different shapes, not just different sizes, even though the topology of the tiling does not change. We call such a tiling “quasi-scalable”.

\subsection{Tiling a large sphere: discovery of the problem}

For “sufficiently large” spheres (i.e., of radius larger than some specified value), we are not aware of any literature prior to 2019. Indeed, we suspect that there has been a widespread assumption that the Isbell tiling could be straightforwardly adapted to tile a sphere of any radius above some threshold value. However, during the “Polymath 16" project, set up to explore the Hadwiger-Nelson problem and related topics, it was noticed that this is not so.

\begin{figure}[!b]
\centering
\includegraphics[scale=0.3]{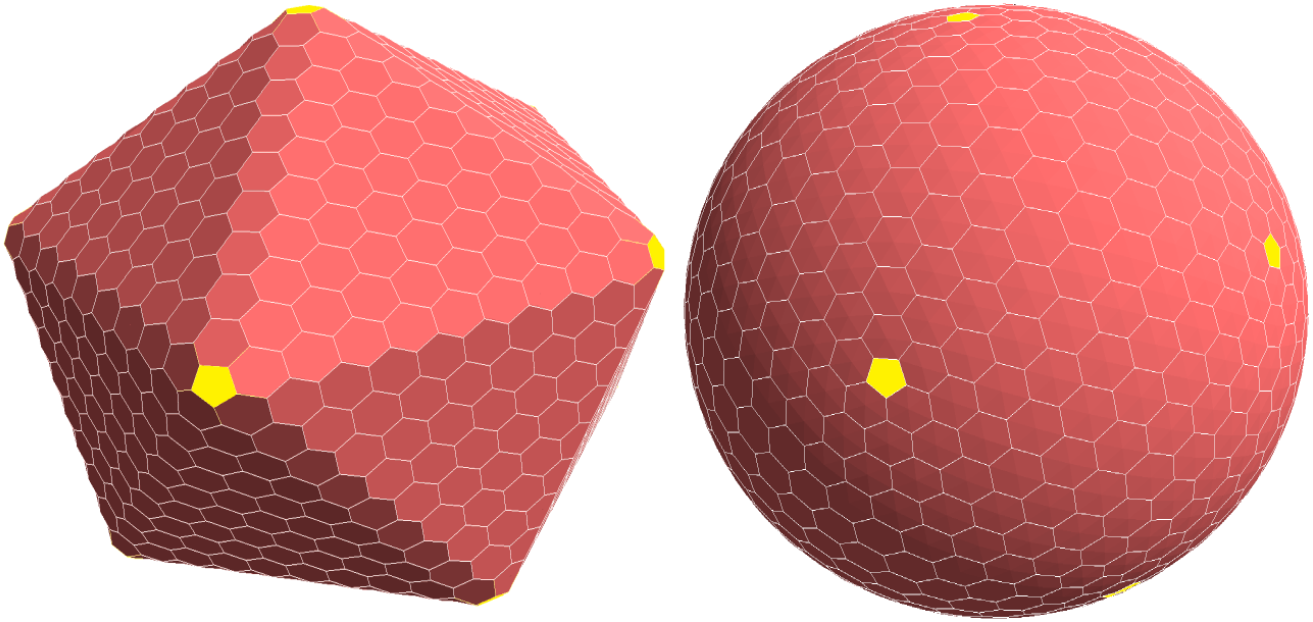}
\caption{Equilateral (left) and spherical (right) Goldberg polyhedra. Image from Wikipedia.}
\label{f4}
\end{figure}

A Goldberg polyhedron is a generalization of the regular truncated icosahedron. Its faces are mostly hexagons, but it has $12$ pentagonal faces. It can be defined as equilateral or spherical (see Figure \ref{f4}); here we shall work mostly with the equilateral form, even though our eventual tilings will be of a sphere. A large Goldberg polyhedron approximates an icosahedron with the pentagons located at its vertices. It will be convenient to have terminology for the faces, edges and vertices of this icosahedron; we shall use the terms icoface, icoedge and icovertex respectively.

One further feature of Goldberg polyhedra must be mentioned. In Figure \ref{f4}, pairs of pentagons are joined by straight-line chains of hexagons. However, one can also create such polyhedra in which the orientation of the hexagonal grids is skewed with respect to the lines between pentagons. The standard convention for naming Goldberg polyhedra incorporates this: $GP(m,n)$ is the Goldberg polyhedron in which a pentagon is reached from another pentagon by stepping $m$ hexagons in one direction, turning $60$ degrees and then stepping $n$ more hexagons.

An attractive way to attempt to $7$-tile the surface of a large sphere, then, is to color a Goldberg polyhedron according to some variant of the Isbell tiling, and then to project the tile boundaries onto the sphere along lines radiating from the origin. Some adjustment is then required in order to compensate for the disparity in distance from the origin of the icovertices and the centers of the icofaces, but that turns out to be simple enough. However, when we come to choose tile colors, the attempt to mimic the Isbell pattern fails. Indeed, recently Ágoston was able to show \cite{ago} that the surface of a sphere of radius greater than $18$ cannot be $7$-colored with a scalable, non-Siamese tiling (i.e. with all the three Thomassen criteria).

In this paper we demonstrate a quasi-scalable, non-Siamese $7$-tiling of a sufficiently large sphere. We note that this is the first case (other than the sphere of radius exactly $\sqrt{\frac{1}{2}}$) of a setting in which scalability definitively matters for chromatic number, i.e. in which the chromatic number when scalability is stipulated is established to be strictly greater than when it is not.

\section{Isbell-tiling a Goldberg polyhedron: reconciling icoedges by deforming tiles}

\subsection{Icoedge mismatches}

\begin{figure}[!b]
\centering
\includegraphics[scale=0.4]{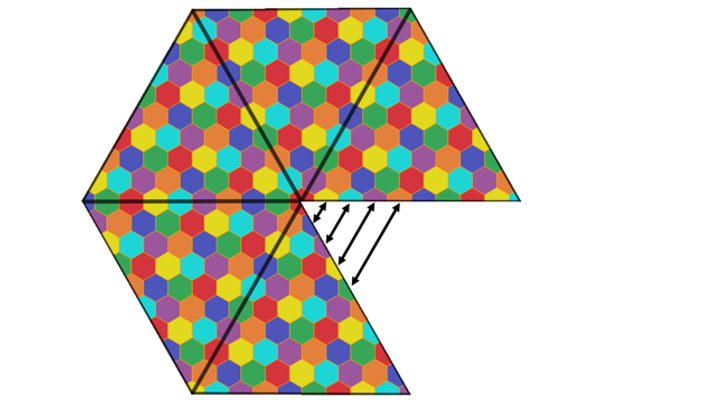}
\caption{The icofaces around a given icovertex cannot all be Isbell-colored: at least one icoedge will not match.}
\label{f5}
\end{figure}

Clearly the tiles of a single icoface of any Goldberg polyhedron, excluding those that straddle the icoedges, can be $7$-colored according to the Isbell pattern. So far so good, but a problem arises when we attempt to reconcile the colorings of icofaces that share an icoedge. The colors of the tiles comprising the first icoface can be transposed arbitrarily, but once we fix our choice we also enforce a unique coloring of the second face. Unfortunately, it turns out that when we continue this for successive icofaces around a given icovertex, the icoedge shared by the first and fifth icofaces does not match up. See Figure \ref{f5}, which uses $GP(9,0)$ for illustration; the same problem arises for all $GP(m,n)$.

\subsection{Reducing icoedge mismatches to a minimal form}

Interestingly, however, this mismatch can be "distributed" across the five icoedges in a manner that makes it quite slight: specifically, at each icoedge just two colors are swapped in one icoface relative to what would match the other icoface.

\begin{figure}[!b]
\centering
\includegraphics[scale=0.4]{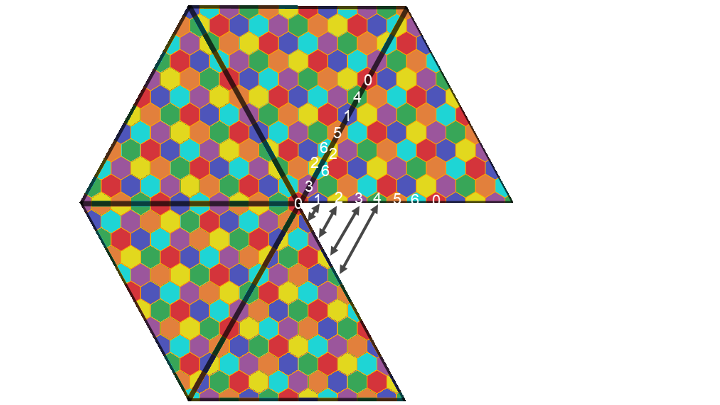}
\caption{In this coloring, all five icoedges are minimally mismatched. (For example, across the horizontal icoedge on the left, orange and yellow are swapped.) As a result, across all five icoedges there are pairs of tiles that are too close together.}
\label{f6}
\end{figure}

To see how this works, see Figure \ref{f6}. Let $swap23$ be the mapping [0,1,2,3,4,5,6] - [0,1,3,2,4,5,6], and let $hex$ be the mapping [0,1,2,3,4,5,6] - [0,3,6,2,5,1,4] that relates one edge of an Isbell-colored icoface to another, going anticlockwise around an icovertex. Then let the sequence of hexagon colors along one side of some icoedge $E$ be [0,1,2,3,4,5,6]. Then five sequential applications of $swap23(hex(E))$ give [0,3,2,6,5,1,4], [0,6,2,4,1,3,5], [0,4,2,5,3,6,1], [0,5,2,1,6,4,3] and finally [0,1,2,3,4,5,6], as required.

Further, note that Figures \ref{f5} and \ref{f6} show Goldberg polyhedra whose second parameter is zero, but this is not required in order for the sequence of five applications of $hex$ and $swap23$ to work. Visually: if we start with the coloring shown in Figure 6, and we slowly rotate the entire honeycomb around the central icovertex, switching colors of hexagons that move across an icoedge according to the Isbell coloring of their new icoface, we will successively generate colorings of (a region of) all $GP(m,n)$. The simplest tiling we have discovered (see later sections) exploits this fact.

\subsection{Resolving a minimal mismatch by introducing unscalable tiles}

\begin{figure}[!b]
\centering
\includegraphics[scale=0.4]{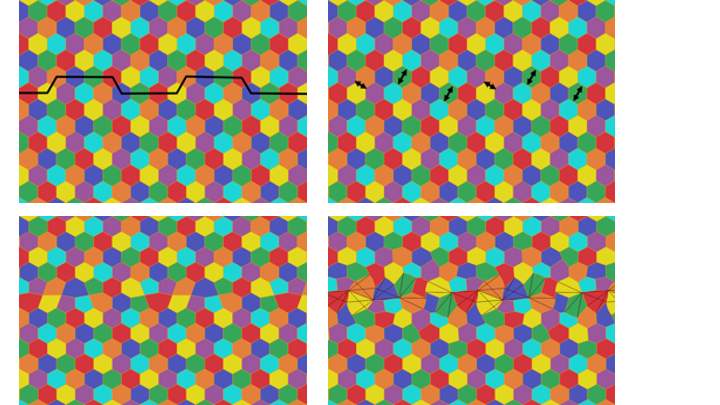}
\caption{Resolution of a minimal mismatch along an icoedge of $GP(m,0)$ using a repeating pattern of unscalable tiles. Top left: purple and cyan tiles are switched on either side of the icoedge. Top right: as a result, a repeating sequence of two pairs of cyan tiles and one pair of purple tiles are too close together. Bottom left: certain inter-tile boundaries along the icoedge are shrunk to a point. Now, many tiles have two same-colored neighbors, but in all such cases one of the contacts is only a single point. Bottom right: various tiles are reshaped to make no tile’s diameter exceed 1 and no same-colored tiles closer than 1. Thin black lines within a tile are constrained to be of unit length.}
\label{f7}
\end{figure}

However, even a minimal mismatch is still a mismatch, giving instances of same-colored tiles less than 1 apart across an icoedge. But it turns out that this can be resolved using the same sort of tactic shown in Figures \ref{f2} and \ref{f3}, with certain inter-tile boundaries being drawn as unit-radius arcs. We proceed in two steps. First we contract certain tile boundaries to points where four tiles meet (see Figure \ref{f7}, bottom left panel). At this point, even though many tiles are either too close or too large, we have eliminated all cases of a tile having two same-colored neighbors with both of which the contact has positive measure. As a result, the coloring can be made valid by moving and curving certain inter-tile boundaries (see Figure \ref{f7}, bottom right panel). Note that we are ignoring, for now, the fact that the tiles shown straddle an icoedge and thus are not all in the same plane, because in the end we will be projecting our construction onto a large sphere, which is locally close to being a plane anywhere. The pattern of deformations repeats every seven hexagons; we call a single repeat a “stitch”.

\begin{figure}[!b]
\centering
\includegraphics[scale=0.3]{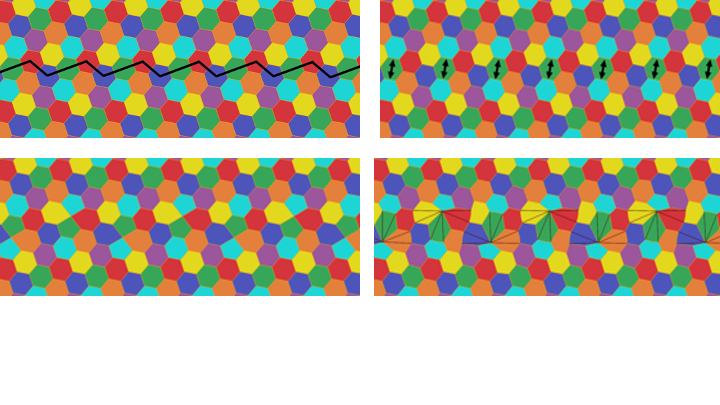}
\caption{Resolution of a minimal mismatch along an icoedge of $GP(4m,2m+1)$ using a repeating pattern of unscalable tiles. Note that the situation is simpler than in Figure \ref{f7}: in the initial tiling (top right panel) there is only one type of over-close pair of tiles, and in the resolved tiling (bottom right), considerably fewer tiles are unscalable.}
\label{f8}
\end{figure}

Our initial work focused on the construction just described, and we were able to proceed all the way to a valid $7$-tiling of a sufficiently large sphere. However, we omit the details, because we later discovered that greater simplicity arises if we work not with $GP(m,0)$ but with $GP(4m,2m+1)$. The corresponding resolution of an icoedge segment is shown in Figure \ref{f8}.

\subsection{Extending the tiling to all 30 icoedges}

In both of the icoedge tilings discussed above, a single stitch possesses rotational symmetry in respect of tile shapes (though not tile colors): in Figure \ref{f7} (bottom right), rotation of a stitch by 180 degrees around the midpoint of the yellow/orange tile boundary superimposes tiles precisely, and the same goes for Figure \ref{f8} with the midpoint of the red/blue boundary. We now impose a new constraint: that the icoedge as a whole must have such rotational symmetry. In other words, the midpoint of an icoedge must coincide with either the midpoint of a stitch or the midpoint between two stitches.

The reason for introducing this constraint is that it ensures that all icovertices are the same: if we can tile a single icovertex, we can tile them all. Moreover, the same applies to colors of the tiles: in these tilings, the rotation by 180 degrees is equivalent to swapping certain pairs of colors, so all we shall need to do is arrange that the colors associated with the icovertices at either end of an icoedge are a pair that is swapped by that icoedge's rotation. Associating each of six colors with two antipodal icovertices always suffices to allow this in a manner compatible with the arrangement of icoedge colorings already described (see Figure \ref{f6}).

\subsection{Resolving the tiling of the icovertices}

\begin{figure}[!b]
\centering
\includegraphics[scale=0.4]{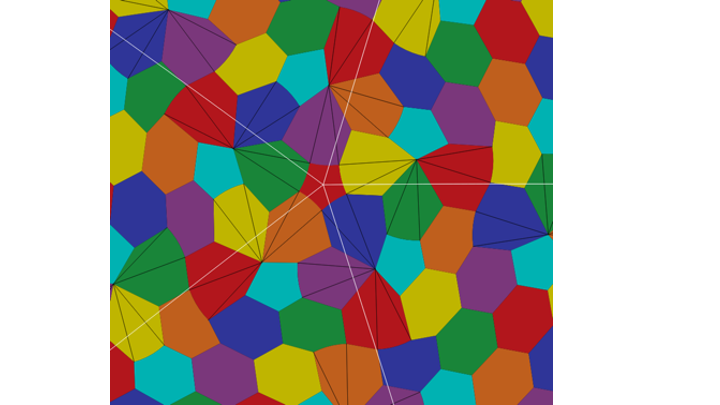}
\caption{A valid 7-tiling of an icovertex region. Observe that the convergence of the icoedges at the central pentagon does not create any new cases of too-close tiles, even though the pentagon is quite small.}
\label{f9}
\end{figure}

In addition to supporting the simplest deformation that gives icoedges a valid coloring, the tiling shown in Figure \ref{f8} is particularly permissive at the icovertices. As shown in Figure \ref{f9}, the stitches at the ends of the five icoedges at an icovertex dovetail naturally in a manner that does not create any new too-close, same-colored tile pairs. This contrasts with our findings for other stitch designs, in which the icovertex regions required more elaborate handling.

Thus, at this point we have a construction for coloring any $GP(4m,2m+1)$, using only seven colors, by deforming certain tiles along its icoedges and at its icovertices, in such a way that any small patch of the polyhedron, if projected onto a plane, is a valid coloring (i.e. contains no pair of points at unit distance).

\section{Projecting the deformed Goldberg polyhedron onto a sphere}

Despite this, it is by no means obvious that our deformed Goldberg polyhedron can be projected onto a sphere in a manner that provides a valid 7-tiling of the sphere. However, as depicted in Figure \ref{f4} (right panel), a Goldberg polyhedron can be defined in a spherical form, with all tile vertices lying on a sphere. This gives cause for optimism, which turns out to be justified.

\subsection{Length of a stitch}

\begin{figure}[!b]
\centering
\includegraphics[scale=0.35]{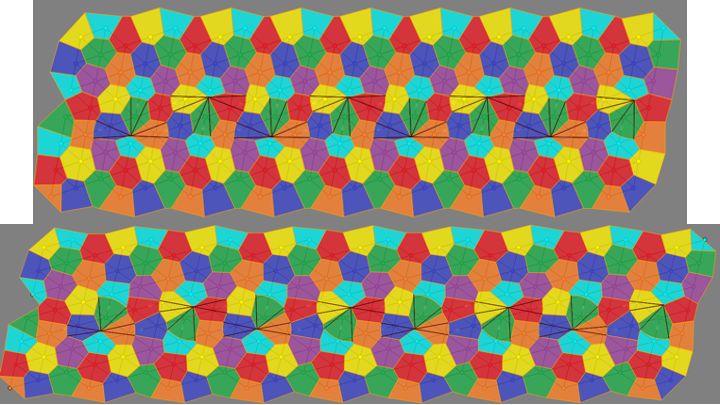}
\caption{A more “compressed” (top) and a more “stretched” (bottom) version of the same icoedge segment.}
\label{f10}
\end{figure}

The tiling shown in Figure \ref{f8} (bottom right panel) contains a repeating sequence of six unscalable tiles, but these sextets are separated by scalable tiles. Thus, we may expect that such a tiling will be quasi-scalable. This is indeed the case: it turns out to be straightforward to construct “long” and “short” versions of the icoedge region differing in stitch length by a factor of more than 1.2 (see Figure \ref{f10}). Clearly all intermediate lengths are therefore also possible.

Thus, if we were only trying to tile the 30 icoedges on a sphere, we would be done: for each value of $m$, there is a range of radii $[r,1.2r]$ for which the tiling corresponding to $GP(4m,2m+1)$ can be drawn on the sphere and $7$-colored. The convexity of the sphere is easily seen not to be able to bring tiles closer together along the icoedge than they are in the plane, even when their boundaries are curved. Thus, since for sufficiently large $m$ we have $m+1 \leq 1.2m$, there is a threshold value of $m$ above which the ranges of radius tilable according to $GP(4m,2m+1)$ for consecutive values of $m$ overlap, so we can $7$-color any sphere larger than that.

\subsection{Filling in the icofaces}

\begin{figure}[!b]
\centering
\includegraphics[scale=0.35]{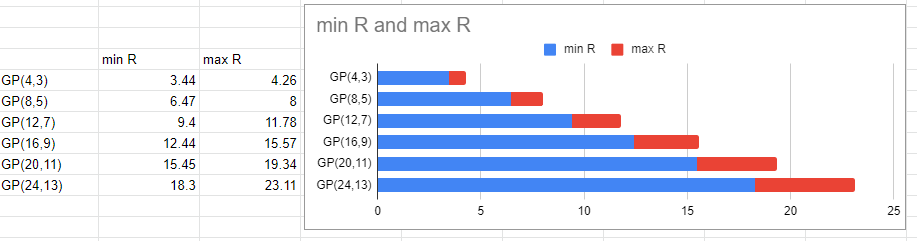}
\caption{The sphere radii that can be $7$-tiled according to our deformation of $GP(4m,2m+1)$, for the first few values of $m$.}
\label{f11}
\end{figure}

It is challenging to define an exact formula for the locations of the tile boundaries of a spherical Goldberg polyhedron, but intuitively it seems unlikely that there would be any difficulty choosing them so as to retain the validity of the Isbell coloring, since the icofaces do not contain any unscalable tiles. To check this, we used iterative numerical methods to  optimize the tiling for $m=1..6$, testing how large and small the sphere radius could be and still allow the tiling to be valid. Our findings are summarized in Figure \ref{f11}, and suggest that any sphere of radius at least $12.44$ can be $7$-colored. The ratio of the maximum and minimum radii for a given value of $m$ appears to remain roughly constant, so we claim that these ranges for $m$ and $m+1$ always overlap for larger $m$. Moreover, as a further sanity check, we ran this same optimization with the additional constraint that any scalable distance between tile vertices is prohibited from lying in the range $[0.99,1.01]$, and this only narrowed the ranges by about 3 percent.

\begin{figure}[!b]
\centering
\includegraphics[scale=0.35]{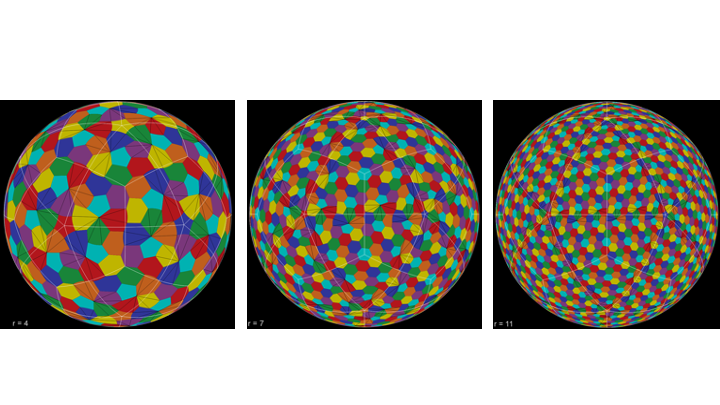}
\caption{$7$-tilings of spheres of radius $4$, $7$ and $11$.}
\label{f12}
\end{figure}

This completes the justification of our main claim, that the surface of any sufficiently large sphere can be $7$-tiled. We conclude (Figure \ref{f12}) with pictures of the whole sphere for $m=1,2$ and $3$, for representative radii.

\section{Acknowledgments}
We thank Aubrey de Grey for extensive editorial assistance, and the participants in the Polymath 16 project for inspiring this work.

\begin{thebibliography}{100}
\bibitem{soi}
A. Soifer. The Mathematical Coloring Book, Springer-Verlag, ISBN 978-0-387-74640-1 (2008).
\bibitem{ago}
P. Ágoston. A lower bound on the number of colours needed to nicely colour a sphere. Proceedings of CCCG 2020, Saskatoon, Canada, August 5–7, 2020.
\bibitem{sim}
G. Malen. Measurable colorings of $\mathbb{S}^2_r$. Geombinatorics 24: 172-180 (2015).
\bibitem{thom}
C. Thomassen. On the Nelson unit-distance coloring problem. Amer. Math. Monthly 106:850-853 (1999).
\bibitem{tow}
S. P. Townsend. Coloring the plane with no monochrome unit. Geombinatorics 14:181-193 (2005).
\bibitem{aphak}
K. Appel, W. Haken. Solution of the Four Color Map Problem. Scientific American 237:108–121 (1977).
\end{thebibliography}
\end{document}